\documentclass[11pt, reqno]{article}
\usepackage{mathrsfs}
\usepackage{amsfonts}
\usepackage{amsmath}
\usepackage{ctex}
\usepackage{CJK}
\usepackage{indentfirst}
\usepackage{cases}
\usepackage{amssymb,amsmath,amscd}
\newcommand{\be}{\begin{equation}}
\newcommand{\ee}{\end{equation}}
\newcommand{\bea}{\begin{eqnarray}}
\newcommand{\eea}{\end{eqnarray}}
\newcommand{\bna}{\begin{eqnarray*}}
\newcommand{\ena}{\end{eqnarray*}}

\begin{document}
\openup 1\jot
\baselineskip=0.5cm

\centerline{\Large Averages of shifted convolution sums for $GL(3) \times GL(2)$
}

\bigskip

\centerline{\large Qingfeng Sun}

\bigskip

\noindent{\small {\bf Abstract} \,
Let $A_f(1,n)$ be the normalized Fourier coefficients of
a $GL(3)$ Maass cusp form $f$ and let
$a_g(n)$ be the normalized Fourier coefficients of a $GL(2)$ cusp form $g$.
Let $\lambda(n)$ be either $A_f(1,n)$ or the triple divisor
function $d_3(n)$.
It is proved that for any $\epsilon>0$, any integer $r\geq 1$ and
$r^{5/2}X^{1/4+7\delta/2}\leq H\leq X$ with
$\delta>0$,
$$
\frac{1}{H}\sum_{h\geq 1}W\left(\frac{h}{H}\right)
\sum_{n\geq 1}\lambda(n)a_g(rn+h)V\left(\frac{n}{X}\right)\ll X^{1-\delta+\epsilon},
$$
where $V$ and $W$ are smooth
compactly supported functions, and the implied constants depend only on
the associated forms and $\epsilon$.
}

\medskip


\medskip

\bigskip

\section{Introduction}
\setcounter{equation}{0}

The shifted convolution sum problems have a long history in analytic number theory.
Nontrivial bounds of various shifted convolution sums have been playing
important roles in many central problems, such as quantum unique ergodicity,
subconvexity and power moments of $L$-functions (see for example \cite{Blomer1}, \cite{DFI1}, \cite{HM}, \cite{HS},
\cite{LLY}, \cite{LS}, \cite{S}).
The first shifted convolution sum involving $GL(3)$ Fourier coefficients was studied
in \cite {P} by Pitt who considered the shifted convolution sum of $d_3(n)$ with the
Fourier coefficients $a_g(n)$ of a holomorphic cusp form $g$, where
$d_3(n)=\sum\limits_{l_1l_2l_3=n \atop l_j\in \mathbb{N}, j=1,2,3}1$ is
the triple divisor function which is the $n$-th coefficient of the cube of the Riemann zeta
function $\zeta^3(s)$.
Recently, Munshi
\cite{munshi1} studied the general $GL(3) \times GL(2)$ shifted convolution sum
\bna
\mathcal {D}_h(X)=\sum_{n\geq 1}A_f(1,n)a_g(n+h)V\left(\frac{n}{X}\right),
\ena
where $A_f(1,n)$ are the Fourier coefficients of a $GL(3)$ Maass
cusp form $f$, $a_g(n)$ are those of a $GL(2)$ Maass or holomorphic cusp
form $g$, $1\leq h\leq X^{1+\epsilon}$ an integer, and $V$ is a smooth
compactly supported function, and succeeded in showing that
\bna
\mathcal {D}_h(X)\ll_{f,g,\epsilon}X^{1-\frac{1}{20}+\epsilon}
\ena
by using the idea of factorizable moduli with the circle method of Jutila's version.
As Munshi remarked in his paper,
``it is expected that extra cancellation can be obtained by averaging over $h$",
which will be the main concern of this paper. In fact, we shall consider the following
averages of $GL(3)\times GL(2)$ shifted convolution sums
\bea
\mathcal {S}(H,X)=\frac{1}{H}\sum_{h\geq 1}W\left(\frac{h}{H}\right)
\sum_{n\geq 1}\lambda(n)a_g(rn+h)V\left(\frac{n}{X}\right),
\eea
where $W$ is another smooth compactly supported function,
$r\geq 1$ is an integer, $\lambda(n)$ is either $A_f(1,n)$ or $d_3(n)$.
Here $f$ is a Maass cusp form for $SL(3,\mathbb{Z})$ and $g$ is a
Maass or holomorphic cusp form for $SL(2,\mathbb{Z})$. Our main result is the following theorem.

\medskip

\noindent {\bf Theorem 1}\, {\it
For any $\epsilon>0$, any integer
$r\geq 1$ and $(rX)^{1/2+\epsilon}\leq H\leq X$, we have
$$\mathcal {S}(H,X)\ll X^{-A}$$ for any $A>0$.
For any $\epsilon>0$, any integer
$r\geq 1$ and $r^{5/2}X^{1/4+7\delta/2}\leq H\leq (rX)^{1/2+\epsilon}$ with
$\delta>0$, we have
\bna
\mathcal {S}(H,X)
\ll X^{1-\delta+\epsilon}.
\ena
Here the implied constants depend only on
the associated forms and $\epsilon$.
}

Recently, averages of shifted convolution sums for $GL(2)$ cusp forms have
been studied in \cite{Blomer},\cite{Lin} and \cite{Singh}.
We note that for the shifted convolution sum in (1.1) without averaging and $\lambda(n)=d_3(n)$,
Munshi's approach for $\mathcal {D}_h(X)$ can also be applied
(see \cite{munshi3}). Moreover, since $d_3(n)\ll n^{\epsilon}$ for any $\epsilon>0$,
we can remove the smooth weight $V$ in (1.1).

\medskip

\noindent {\bf Theorem 2}\, {\it Assume that $a_g(n)\ll n^{\theta+\epsilon}$ for any $\epsilon>0$.
For any $\epsilon>0$, any integer
$r\geq 1$ and $(rX)^{1/2+\epsilon}\leq H\leq X$, we have
$$
\frac{1}{H}\sum_{h\geq 1}W\left(\frac{h}{H}\right)
\sum_{n\leq X}d_3(n)a_g(rn+h)\ll X^{-A}
$$ for any $A>0$.
For any $\epsilon>0$, any integer $r\geq 1$ and
$r^{5/2}X^{1/4+6\delta}(rX)^{5\theta/2}\leq H\leq (rX)^{1/2+\epsilon}$ with
$\delta>0$, we have
\bna
\frac{1}{H}\sum_{h\geq 1}W\left(\frac{h}{H}\right)
\sum_{n\leq X}d_3(n)a_g(rn+h)
\ll_{g,\epsilon} X^{1-\delta+\epsilon}.
\ena
}

Note that we can take $\theta=0$ for $g$ a holomorphic cusp form
and $\theta=7/64$ for $g$ a Maass cusp
form (see \cite{K}).

\medskip

\section{The circle method and Voronoi formulas}
\setcounter{equation}{0}
\medskip

\noindent {\bf 2.1 The circle method}
\medskip

As usual, denote $\delta(n)=\left\{\begin{array}{ll}
1,&\mbox{if $n=0$},\\
0,&\mbox{otherwise}.
\end{array}\right.$

\medskip

\noindent {\bf Lemma 1 (\cite{H})}\, {\it For any $P>1$ there is a positive constant
$c_P$, and a smooth function $h(x,y)$ defined on $(0,\infty)\times \mathbb{R}$,
such that
$$
\delta(n)=\frac{c_P}{P^2}\sum_{q=1}^{\infty}\;\sideset{}{^*}\sum_{c \bmod q}
e\left(\frac{cn}{q}\right)h\left(\frac{q}{P},\frac{n}{P^2}\right)
$$
for $n\in \mathbb{Z}$. Here the $*$ over the sum indicates that $c$ and $q$
are coprime. The constant $c_P=1+O_A(P^{-A})$ for any $A>0$. Moreover,
$h(x,y)\ll x^{-1}$ for all $y$, and $h(x,y)$ is nonzero only when
$x\leq \max\{1, 2|y|\}$. The smooth function $h(x,y)$ satisfies
\bea
x^i \frac{\partial^i h}{\partial x^i}(x,y)\ll_i x^{-1} \quad \mbox{and} \quad
\frac{\partial h}{\partial y}(x,y)=0
\eea
for $x\leq 1$ and $|y|\leq x/2$. And also for $|y|\geq x/2$, we have
\bea
x^iy^j \frac{\partial^{i+j} h}{\partial x^i\partial y^j}(x,y)\ll_{i,j} x^{-1}.
\eea
}

We will apply Lemma 1 for larger $H$ using the fact that we can choose $P=\sqrt{Y}$
to detect the equation $n=0$ for integers in the range $|n|\leq Y$.
For small $H$, Lemma 1 is not efficient to obtain savings (for small $q$) in our problem and we will
apply Jutila's variation of the circle method (\cite{J1})
which gives an approximation for
$
I_{[0,1]}(x)=\left\{\begin{array}{ll}
1,&x\in [0,1],\\
0,&\mbox{otherwise},
\end{array}
\right.
$
where
$
I_{S}(x)
$
is the characteristic function of the set $S$. We have the following result
(for a proof see \cite{munshi1}, Lemma 4).

\noindent {\bf Lemma 2}\, {\it
Let $\mathfrak{Q}\subset [1,Q]$, $Q>0$ and $Q^{-2}\leq \eta \leq Q^{-1}$. Define
\begin{equation}
\widetilde{I}_{\mathfrak{Q},\eta}(x)=\frac{1}{2\eta L}
\sum_{q\in \mathfrak{Q}}\,\sideset{}{^*}\sum_{c \bmod q}
I_{\left[\frac{c}{q}-\eta,\frac{c}{q}+\eta\right]}(x),
\end{equation}
where $L=\sum_{q\in \mathfrak{Q}}\phi(q)$.
Then for any $\epsilon>0$,
\begin{equation}
\int\limits_0^1\left|1-\widetilde{I}_{\mathfrak{Q},\eta}(\beta)\right|^2\mathrm{d}\beta\ll
\frac{Q^{2+\epsilon}}{\eta L^2}.
\end{equation}
}

\bigskip

\noindent{\bf 2.2 $GL(2)$ Voronoi formulas}
\medskip

For notational simplicity, we assume that $g$ is a Hecke-Maass cusp form
for $SL(2,\mathbb{Z})$ with Laplace eigenvalue
$1/4+\mu^2$ and normalized Fourier coefficients $a_g(m)$.

\noindent {\bf Lemma 3 (\cite{M})}\, {\it Let $\psi(y)\in C_c^\infty(0,\infty)$. For $(c,q)=1$, we have
\bna
\sum_{m=1}^{\infty}a_g(m)e\left(\frac{cm}{q}\right)\psi(m)=\frac{1}{q}\sum_{\pm}\sum_{m=1}^{\infty}
a_g(\mp m)e\left(\pm \frac{\overline{c} m}{q}\right)
\Psi^{\pm}\left(\frac{m}{q^2}\right),
\ena
where $\overline{c}$ denote the multiplicative inverse of $c\bmod q$, and
\bea
\Psi^{-}\left(y\right)&=&-\frac{\pi}{\cosh(\pi \mu)}
\int_0^{\infty}\psi(v)(Y_{2i \mu}+Y_{-2i\mu})(4\pi\sqrt{y v})\mathrm{d}v,\\
\Psi^{+}\left(y\right)&=&4\cosh(\pi \mu)
\int_0^{\infty}\psi(v)K_{2i\mu}(4\pi\sqrt{y v})\mathrm{d}v.
\eea
}

If $\psi(y)$ is a smooth function of
compact support in $[AY,BY]$,
where $Y>0$ and $B>A>0$,
satisfying $\psi^{(j)}(y)\ll_{A,B,j}Y^{-j} $ for any integer $j\geq 0$, then
for any fixed $\epsilon>0$ and $yY\gg Y^{\epsilon}$,
$\Psi^{\pm}(y)$ are negligibly small. For $yY\ll Y^{\epsilon}$, we
have the trivial bound
$
\Psi^{\pm}(y)\ll_{g,\epsilon} Y^{1+\epsilon}.
$

\medskip

\noindent{\bf 2.3 $GL(3)$ Voronoi formulas}

\medskip

Let $f$ be a Hecke-Maass cusp form of type $(\nu_1,\nu_2)$ for
$SL(3,\mathbb{Z})$ with normalized Fourier coefficients $A_f(n_1,n_2)$.
Denote
$
\mu_1=-\nu_1-2\nu_2+1, \mu_2=-\nu_1+\nu_2, \mu_3=\nu_1+\nu_2-1.
$
The generalized Ramanujan conjecture asserts that $\mathrm{Re}(\mu_j)=0$, $1\leq j\leq 3$,
while the current record bound due to
Luo, Rudnick and Sarnak \cite{LRS} is
$
|\mathrm{Re}(\mu_j)| \leq \frac{1}{2}-\frac{1}{10},  1\leq j\leq 3.
$

Let $\varphi(y)$ be a smooth function compactly supported on $(0,\infty)$ and
denote by $\widetilde{\varphi}(s)$
the Mellin transform of $\varphi(y)$.
For $k=0,1$, we define
\bea
\Phi_k(y):=\int\limits_{\mathrm{Re}(s)=\sigma}(\pi^3 y)^{-s}\prod_{j=1}^3
\frac{\Gamma\left(\frac{1+s+\mu_j+2k}{2}\right)}
{\Gamma\left(\frac{-s-\mu_j}{2}\right)}
\widetilde{\varphi}(-s-k)\mathrm{d}s
\eea
with $\sigma>\max\limits_{1\leq j\leq 3}\{-1-\mathrm{Re}(\mu_j)-2k\}$. Set
\bea
\Phi^{\pm}(y)=\Phi_0(y)\pm\frac{1}{i\pi^3y}\Phi_1(y).
\eea
Then we have the following Voronoi formula.

\noindent {\bf Lemma 4 (\cite{GL1}, \cite{MS})} \,{\it
Let $\varphi(y) \in C_c^\infty(0,\infty)$.
For $(c,q)=1$ we have
\bna
\sum_{n\geq 1}A_f(1,n)e\left(\frac{cn}{q}\right)\varphi(n)=\frac{q\pi^{-\frac{5}{2}}}{4i}
\sum_{\pm}\sum_{n_1|q}
\sum_{n_2=1}^{\infty}
\frac{A_f(n_2,n_1)}{n_1n_2}S\left(\overline{c},\pm n_2;\frac{q}{n_1}\right)\Phi^{\pm}\left(\frac{n_1^2n_2}{q^3}\right),
\ena
where $\overline{c}$ denote the multiplicative inverse of $c\bmod q$ and $S(m,n;c)$ is the classical Kloosterman sum.
}

Next we state the Voronoi formula for $d_3(n)$ in Li's version (see \cite{Li3}). Set
$
\sigma_{0,0}(k,l)=\sum\limits_{d_1|l \atop d_1>0}\sum\limits_{d_2|\frac{l}{d_1} \atop d_2>0, (d_2,k)=1}1.
$
Let $\gamma:=\lim\limits_{s\rightarrow 1}\left(\zeta(s)-\frac{1}{s-1}\right)$
be the Euler constant and $\gamma_1:=-\frac{\mathrm{d}}{\mathrm{d}s}\left.\left(\zeta(s)-\frac{1}{s-1}\right)\right|_{s=1}$
be the Stieltjes constant. For $\omega(y)\in C_c(0,\infty)$, $k=0,1$ and $\sigma>-1-2k$, set
\bea
\Omega_k(y)=\frac{1}{2\pi i}\int\limits_{\mathrm{Re}(s)=\sigma}\left(\pi^3y\right)^{-s}
\frac{\Gamma\left(\frac{1+s+2k}{2}\right)^3}{\Gamma\left(\frac{-s}{2}\right)^3}
\widetilde{\omega}(-s-k)\mathrm{d}s
\eea
with $\widetilde{\omega}(s)=\int_0^{\infty}\omega(u)u^{s-1}\mathrm{d}u$
the Mellin transform of $\omega$, and
\bea
\Omega^{\pm}(y)=\Omega_0(y)\pm \frac{1}{i\pi^3y}\Omega_1(y).
\eea

\medskip

\noindent{\bf Lemma 5}\, {\it Let $\omega(y)\in C_c^{\infty}(0,\infty)$.
For $(c,q)=1$ and $c\overline{c}\equiv 1(\bmod q)$ we have
\bna
&&\sum_{n\geq 1}d_3(n)e\left(\frac{c n}{q}\right)\omega(n)\\
&=&
\frac{q}{2\pi^{\frac{3}{2}}}\sum_{\pm}\sum_{n|q}\sum_{m\geq 1}\frac{1}{nm}
\sum_{n_1|n}\sum_{n_2|\frac{n}{n_1}}\sigma_{0,0}\left(\frac{n}{n_1n_2},m\right)
S\left(\pm m,\overline{c};\frac{q}{n}\right)\Omega^{\pm}\left(\frac{mn^2}{q^3}\right)\\
&&+\frac{1}{2q^2}\widetilde{\omega}(1)\sum_{n|q}n\tau(n)P_2(n,q)S\left(0,\overline{c};\frac{q}{n}\right)
\\
&&+\frac{1}{2q^2}\widetilde{\omega}'(1)\sum_{n|q}n\tau(n)P_1(n,q)S\left(0,\overline{c};\frac{q}{n}\right)
\\
&&+\frac{1}{4q^2}\widetilde{\omega}''(1)\sum_{n|q}n\tau(n)S\left(0,\overline{c};\frac{q}{n}\right),
\ena
where
$
P_1(n,q)=\frac{5}{3}\log n-3\log q+3\gamma-\frac{1}{3\tau(n)}\sum_{d|n}\log d,
$
and
\bna
P_2(n,q)&=&\left(\log n\right)^2-5\log q \log n
+\frac{9}{2}(\log q)^2+3\gamma^2-3\gamma_1+7\gamma \log n-9\gamma \log q\nonumber\\
&&+\frac{1}{\tau(n)}\left(\left(\log n+\log q-5\gamma\right)\sum_{d|n}\log d-\frac{3}{2}
\sum_{d|n}(\log d)^2\right).
\ena
}

The functions $\Phi^{\pm}(y)$ (also $\Omega^{\pm}(y)$) have the following properties
(see Sun \cite{Sun} for proof).

\noindent{\bf Lemma 6} \,{\it
Suppose that $\varphi(y)$ is a smooth function of
compact support in $[AY,BY]$,
where $Y>0$ and $B>A>0$,
satisfying $\varphi^{(j)}(y)\ll_{A,B,j}P^j $ for any integer $j\geq 0$.
Then for $y>0$ and any integer $\ell\geq 0$, we have
$$
\Phi^{\pm}(y)\ll_{A,B,\ell,\epsilon}(yY)^{-\epsilon}(PY)^3\left(\frac{y}{P^3Y^2}\right)^{-\ell}.
$$
}

By Lemma 6, for any fixed $\epsilon>0$ and $yY\gg Y^{\epsilon}(P Y)^3$,
$\Phi^{\pm}(y)$ are negligibly small. For $yY\ll Y^{\epsilon}(P Y)^3$, we can
shift the contour of integration in (2.7) to $\sigma=-3/5+\epsilon$ with $\epsilon>0$ to get
\bea
\Phi^{\pm}(y) &\ll& (y Y)^{\frac{3}{5}-\epsilon}P Y.
\eea

\section{Proof of Theorem 1}
\setcounter{equation}{0}
\medskip

We write
\bea
\mathcal {S}(H,X)=\frac{1}{H}\sum_{h\geq 1}W\left(\frac{h}{H}\right)
\sum_{n\geq 1}\lambda(n)V\left(\frac{n}{X}\right)\sum_{m\geq 1}a_g(m)
\phi\left(\frac{m}{rX+h}\right)\delta\left(rn+h-m\right),
\eea
where $\phi(y)$ is a smooth function compactly supported in $[1/2,5/2]$, which equals 1
on $[1,2]$ and satisfies $\phi^{(j)}(y)\ll_j 1$. Taking $P=\sqrt{6rX}$ and
applying the circle method in Lemma 1, we have
\bna
\mathcal {S}(H,X)
&=&\frac{c_P}{HP^2}\sum_{q\leq P}\;\sideset{}{^*}\sum_{c \bmod q}
\sum_{n\geq 1}\lambda(n)e\left(\frac{c rn}{q}\right)V\left(\frac{n}{X}\right)
\sum_{m\geq 1}a_g(m)e\left(-\frac{cm}{q}\right)\nonumber\\
&&\sum_{h\geq 1}e\left(\frac{ch}{q}\right)W\left(\frac{h}{H}\right)\phi\left(\frac{m}{rX+h}\right)
h\left(\frac{q}{P},\frac{rn+h-m}{P^2}\right).
\ena
Applying Poisson summation to the $h$-sum (see Theorem 4.4 in \cite{IK}), we have
\bna
\mbox{$h$-sum}&=&\sum_{\gamma\bmod q}e\left(\frac{c \gamma}{q}\right)
\sum_{h\equiv \gamma \bmod q} W\left(\frac{h}{H}\right)
\phi\left(\frac{m}{rX+h}\right)h\left(\frac{q}{P},\frac{rn+h-m}{P^2}\right)\\
&=&\frac{1}{q}\sum_{\gamma\bmod q}e\left(\frac{c \gamma}{q}\right)
\sum_{h\in \mathbb{Z}}e\left(\frac{h \gamma}{q}\right)
\int_{\mathbb{R}}W\left(\frac{x}{H}\right)
\phi\left(\frac{m}{rX+x}\right)h\left(\frac{q}{P},\frac{rn+x-m}{P^2}\right)e\left(-\frac{h x}{q}\right)\mathrm{d}x\\
&=&H\sum_{h\in \mathbb{Z}\atop h\equiv -c\bmod q}\mathcal {I}(h,n,m,q),
\ena
where
\bna
\mathcal {I}(h,n,m,q)=\int_{\mathbb{R}}W\left(x\right)
\phi\left(\frac{m}{rX+Hx}\right)h\left(\frac{q}{P},\frac{rn+Hx-m}{P^2}\right)
e\left(-\frac{h H x}{q}\right)\mathrm{d}x.
\ena
Note that the condition $h\equiv -c \bmod q$ implies that $(h,q)=1$. Then for $h=0$, we have
$q=1$. For $h\neq 0$, by partial integration $j$ times and (2.1)-(2.2), we have
\bna
\mathcal {I}(h,n,m,q)
\ll_j\frac{P}{q}\left(\frac{|h| H}{q}\right)^{-j}
\left(1+\frac{H}{rX+H}+\frac{P}{q}\frac{H}{P^2}\right)^j
\ll_j\frac{P}{q}\left(\frac{P}{H|h|}\right)^j.
\ena
Thus the contribution from $|h|\geq P^{1+\epsilon}/H$
is negligible. In particular, if $H>(rX)^{\frac{1}{2}+\epsilon}$, we have
\bna
\mathcal {S}(H,X)&=&\frac{c_P}{P^2}
\sum_{n\geq 1}\lambda(n)V\left(\frac{n}{X}\right)
\sum_{m\geq 1}a_g(m)\mathcal {I}(0,n,m,1)\\
&=& \frac{c_P}{P^2}\int_{\mathbb{R}}W\left(x\right)
\sum_{n\geq 1}\lambda(n)V\left(\frac{n}{X}\right)
\sum_{m\geq 1}a_g(m)
\phi\left(\frac{m}{rX+Hx}\right)
h\left(\frac{1}{P},\frac{rn+Hx-m}{P^2}\right)\mathrm{d}x\\
&\ll&X^{-A}
\ena
for any $A>0$. Here we have used the fact of Booker \cite{B} that for $\pi$
an automorphic representation of $GL_r(\mathbb{A}_\mathbb{Q})$ whose $L$-function
$L(s,\pi)=\sum_{n\geq 1}\lambda_{\pi}(n)n^{-s}$ is entire, and $F$ a Schwartz
function on $(0,\infty)$,
\bea
\sum_{n\geq 1}\lambda_{\pi}(n)F\left(\frac{n}{X}\right)\ll_{\pi,A,F}X^{-A}
\eea
for any $A>0$.

For $H\leq(rX)^{\frac{1}{2}+\epsilon}$, we write (3.1) as
\bna
\mathcal {S}(H,X)&=&\frac{1}{H}\sum_{h\geq 1}W\left(\frac{h}{H}\right)
\sum_{n\geq 1}\lambda(n)V\left(\frac{n}{X}\right)\sum_{m\geq 1}a_g(m)
\phi\left(\frac{m}{rX+h}\right)\int_0^1e\left(\left(rn+h-m\right)\alpha\right)\mathrm{d}\alpha\\
&=&\frac{1}{H}\sum_{h\geq 1}W\left(\frac{h}{H}\right)\int_0^1 e(\alpha h)
\sum_{n\geq 1}\lambda(n)e(\alpha r n)V\left(\frac{n}{X}\right)\sum_{m\geq 1}a_g(m)
e(-\alpha m)\phi\left(\frac{m}{rX+h}\right)\mathrm{d}\alpha.
\ena
By Lemma 2 we shall approximate $\mathcal {S}(H,X)$ by
\bna
\widetilde{\mathcal {S}}(H,X)&=&\frac{1}{H}\sum_{h\geq 1}W\left(\frac{h}{H}\right)\int_0^1
\widetilde{I}_{\mathfrak{Q},\eta}(\alpha)e(\alpha h)
\sum_{n\geq 1}\lambda(n)e(\alpha r n)V\left(\frac{n}{X}\right)\\
&&\sum_{m\geq 1}a_g(m)
e(-\alpha m)\phi\left(\frac{m}{rX+h}\right)\mathrm{d}\alpha,
\ena
where $\widetilde{I}_{\mathfrak{Q},\eta}(\alpha)(x)$ is defined in (2.3).
Then by Cauchy's inequality and (2.4),
\bna
\mathcal {S}(H,X)-\widetilde{\mathcal {S}}(H,X)
&\ll&\frac{1}{H}\sum_{h\geq 1}W\left(\frac{h}{H}\right)\int_0^1
\left|1-\widetilde{I}_{\mathfrak{Q},\eta}(\alpha)\right|\\
&&\left|\sum_{n\geq 1}\lambda(n)e(\alpha rn)V\left(\frac{n}{X}\right)\right|
\left|\sum_{m\geq 1}a_g(m)
e(-\alpha m)
\phi\left(\frac{m}{rX+h}\right)\right|\mathrm{d}\alpha\\
&\ll_{g,\epsilon}& (r X)^{\frac{1}{2}+\epsilon}
\left(\int_0^1
\left|1-\widetilde{I}_{\mathfrak{Q},\eta}(\alpha)\right|^2\mathrm{d}\alpha\right)^{1/2}\\
&&\left(\int_0^1\left|\sum_{n\geq 1}\lambda(n)e(\alpha rn)V\left(\frac{n}{X}\right)\right|
\mathrm{d}\alpha\right)^{1/2}\\
&\ll_{f,g,\epsilon}& (r X)^{\epsilon}r^{1/2}X\left(\frac{Q^{2+\epsilon}}{\eta L^2}\right)^{1/2}\\
&\ll_{f,g,\epsilon}& (r X)^{\epsilon}\frac{r^{1/2}X}{\sqrt{\eta}Q}
\ena
since $L\gg Q^{2-\epsilon}$, where we have used the Rankin-Selberg estimate
$\sum_{n\leq Y}|A_f(1,n)|^2\ll_{f,\epsilon}Y^{1+\epsilon}$ and the uniform bound in
$\alpha\in \mathbb{R}$
\bna
\sum_{m\geq 1}a_g(m)
e(-\alpha m)
\phi\left(\frac{m}{Y}\right)\ll_{g,\epsilon}Y^{1/2+\epsilon}.
\ena
Taking $\eta=(rX+H)^{-1}$ we obtain
\bea
\mathcal {S}(H,X)=\widetilde{\mathcal {S}}(H,X)
+O\left((r X)^{\epsilon}\frac{rX^{\frac{3}{2}}}{Q}\right).
\eea
Then we only need to estimate $\widetilde{\mathcal {S}}(H,X)$.
Changing variable $\alpha\rightarrow\frac{c}{q}+\beta$, we have
\bna
\widetilde{\mathcal {S}}(H,X)
:=\frac{1}{2\eta}\int_{-\eta}^{\eta}\widetilde{\mathcal {S}}_{\beta}(H,X)\mathrm{d}\beta,
\ena
where
\bea
\widetilde{\mathcal {S}}_{\beta}(H,X)&=&
\frac{1}{HL}\sum_{q\in \mathfrak{Q}}\;\sideset{}{^*}\sum_{c \bmod q}
\sum_{n\geq 1}\lambda(n)e\left(\frac{c r n}{q}\right)
V\left(\frac{n}{X}\right)e\left(\beta r n\right)\sum_{m\geq 1}a_g(m)
e\left(-\frac{c m}{q}\right)
e\left(-\beta m\right)\nonumber\\
&&\sum_{h\geq 1}e\left(\frac{c h}{q}\right)W\left(\frac{h}{H}\right)
\phi\left(\frac{m}{rX+h}\right)e\left(\beta h\right).
\eea
Applying Poisson summation to the $h$-sum, we have
\bna
\mbox{$h$-sum}&=&\sum_{\gamma\bmod q}e\left(\frac{c \gamma}{q}\right)
\sum_{h\equiv \gamma \bmod q} W\left(\frac{h}{H}\right)
\phi\left(\frac{m}{rX+h}\right)e\left(\beta h\right)\\
&=&\frac{1}{q}\sum_{\gamma\bmod q}e\left(\frac{c \gamma}{q}\right)
\sum_{h\in \mathbb{Z}}e\left(\frac{h \gamma}{q}\right)
\int_{\mathbb{R}}W\left(\frac{x}{H}\right)
\phi\left(\frac{m}{rX+x}\right)e\left(\beta x\right)e\left(-\frac{h x}{q}\right)\mathrm{d}x\\
&=&H\sum_{h\in \mathbb{Z}\atop h\equiv -c\bmod q}I_{\beta}(h,m,q),
\ena
where
\bna
I_{\beta}(h,m,q)=\int_{\mathbb{R}}W\left(x\right)
\phi\left(\frac{m}{rX+Hx}\right)e\left(\beta H x\right)e\left(-\frac{h H x}{q}\right)\mathrm{d}x.
\ena

Now we choose the set of moduli $\mathfrak{Q}$ as the prime set
\bna
\mathfrak{Q}=\{q: q\in [Q/2,Q] \ \text{is prime and}\, (q,r)=1\}.
\ena
Then the requirement $L\gg_{\epsilon} Q^{2-\epsilon}$ is satisfied and $h\neq 0$
since the condition $h\equiv -c \bmod q$ implies that $(h,q)=1$.
By partial integration $j$ times we have
$
I_{\beta}(h,m,q)\ll_j \left(|h| H/q\right)^{-j}
$
since $|\beta|\leq \eta=(rX+H)^{-1}$.
Thus contribution from $|h|\geq Q^{1+\epsilon}/H$
is negligible and
\bea
\mbox{$h$-sum}=H\sum_{1\leq |h|\leq \frac{Q^{1+\epsilon}}{H}\atop h\equiv -c\bmod q}I_{\beta}(h,m,q)
+O((r X)^{-A})
\eea
for any $A>0$. Plugging (3.5) into (3.4), we need to estimate
\bna
&&\frac{1}{L}\sum_{1\leq |h|\leq \frac{Q^{1+\epsilon}}{H}}
\sum_{q\in \mathfrak{Q}\atop (q,h)=1}
\sum_{n\geq 1}\lambda(n)e\left(-\frac{h r n}{q}\right)
V\left(\frac{n}{X}\right)e\left(\beta rn\right)\sum_{m\geq 1}a_g(m)
e\left(\frac{h m}{q}\right)
e\left(-\beta m\right)I_{\beta}(h,m,q)\\
&:=&\int_{\mathbb{R}}W\left(x\right)
e\left(\beta H x\right)
\widetilde{\mathcal {S}}_{\beta,x}(H,X)\mathrm{d}x,
\ena
where
\bea
\widetilde{\mathcal {S}}_{\beta,x}(H,X)&=&
\frac{1}{L}\sum_{1\leq |h|\leq \frac{Q^{1+\epsilon}}{H}}
\sum_{q\in \mathfrak{Q}\atop (q,h)=1}e\left(-\frac{h H x}{q}\right)
\sum_{n\geq 1}\lambda(n)e\left(-\frac{h rn}{q}\right)
V\left(\frac{n}{X}\right)e\left(\beta r n\right)\nonumber\\
&&\sum_{m\geq 1}a_g(m)
e\left(\frac{h m}{q}\right)\phi\left(\frac{m}{rX+Hx}\right)
e\left(-\beta m\right).
\eea
We first apply $GL(2)$ Voronoi formula in Lemma 3 to the $m$-sum in (3.6) to get
\bea
\mbox{$m$-sum}=\frac{1}{q}\sum_{\pm}\sum_{m=1}^{\infty}
a_g(\mp m)e\left(\pm \frac{\overline{h} m}{q}\right)
\Psi_{\beta,x}^{\pm}\left(\frac{m}{q^2}\right),
\eea
where $\Psi_{\beta,x}^{\pm}\left(y\right)$ are defined in (2.5)-(2.6) with
$\psi(y)=\phi\left(\frac{y}{rX+Hx}\right)
e\left(-\beta y\right)$.
Note that
\bna
\frac{d^j}{dy^j}\left\{\phi\left(\frac{y}{rX+Hx}\right)
e\left(-\beta y\right)\right\}\ll
\left(\frac{1}{rX+H}+|\beta|\right)^j\ll
\left(\frac{1}{rX}\right)^j.
\ena
Thus the contribution from $|m|\gg Q^2(rX)^{\epsilon}/(rX)$
in (3.7) is negligible. For $|m|\ll Q^2(rX)^{\epsilon}/(rX)$,
we have the trivial bound
$\Psi_{\beta,x}^{\pm}\left(\frac{m}{q^2}\right)\ll (rX)^{1+\epsilon}$.

Next we want to apply the
Voronoi formulas to the $n$-sum
in (3.6).

\medskip

{\bf Case (i) $\lambda(n)=A_f(1,n)$.}
We apply the $GL(3)$ Voronoi formula in Lemma 4 to the $n$-sum
in (3.6) to get
\bea
&&\sum_{n\geq 1}A_f(1,n)e\left(-\frac{h r n}{q}\right)
V\left(\frac{n}{X}\right)e\left(\beta r n\right)\nonumber\\&=&
\frac{q\pi^{-\frac{5}{2}}}{4i}\sum_{\pm}\sum_{n_1|q}\sum_{n_2=1}^{\infty}
\frac{A_f(n_2,n_1)}{n_1n_2}S\left(-\overline{hr},\pm n_2;\frac{q}{n_1}\right)
\Phi_{\beta}^{\pm}\left(\frac{n_1^2n_2}{q^3}\right),
\eea
where $\Phi_{\beta}^{\pm}\left(y\right)$ are defined in (2.7)-(2.8) with
$\varphi(y)=V\left(y/X\right)e\left(\beta r y\right)$. Note that
$
\frac{d^j}{dy^j}\left\{V\left(\frac{y}{X}\right)
e\left(\beta r y\right)\right\}\ll X^{-j}
$
for any $j\geq 0$. By Lemma 6, one sees that the contribution from
$n_1^2n_2\gg Q^3X^{\epsilon}/X$ in (3.8) is negligible.
For $n_1^2n_2\ll Q^3X^{\epsilon}/X$, by (2.11) we get
\bea
\Psi_{\beta}^{\pm}\left(\frac{n_1^2n_2}{q^3}\right)\ll_{f,\epsilon}
\left(\frac{Xn_1^2n_2}{q^3}\right)^{\frac{3}{5}-\epsilon}.
\eea
By (3.6)-(3.9) and Weil's bound for Kloosterman sums, we conclude that
\bea
\widetilde{\mathcal {S}}_{\beta,x}(H,X)&\ll_{f,g,\epsilon}&
\frac{1}{L}\sum_{\pm}\sum_{1\leq |h|\leq \frac{Q^{1+\epsilon}}{H}}
\sum_{q\in \mathfrak{Q}}\sum_{|m|\ll Q^2(rX)^{\epsilon}/(rX)}
|a_g(\mp m)|(r X)^{1+\epsilon}\nonumber\\
&&\sum_{n_1|q}\,\sum_{n_2\ll Q^3X^{\epsilon}/(n_1^2X)}
\frac{|A_f(n_2,n_1)|}{n_1n_2}\left(\frac{q}{n_1}\right)^{1/2}
\left(\frac{Xn_1^2n_2}{q^3}\right)^{\frac{3}{5}-\epsilon}\nonumber\\
&\ll_{f,g,\epsilon}&\frac{X^{\frac{3}{5}+\epsilon}Q}{H}\sum_{q\in \mathfrak{Q}}q^{-\frac{13}{10}}
\sum_{n_1|q}n_1^{-\frac{3}{10}}\sum_{n_2\ll Q^3X^{\epsilon}/(n_1^2X)}
|A_f(n_2,n_1)|n_2^{-\frac{2}{5}}\nonumber\\
&\ll_{f,g,\epsilon}&
\frac{(rX)^{\epsilon}Q^{5/2}}{H},
\eea
where we have used the Rankin-Selberg estimates $\sum_{|m|\ll N}
|a_g(\mp m)|\ll_g N$ and
$
\sum_{n_2\leq N}|A_f(n_1,n_2)|\ll_f N|n_1|.
$
Taking $$Q=(rH)^{2/7}X^{3/7}.$$
Then for $\lambda(n)=A_f(1,n)$ Theorem 1 follows from (3.3) and (3.10).

\medskip

{\bf Case (ii) $\lambda(n)=d_3(n)$.}
Applying Lemma 5 to the $n$-sum
in (3.6) we get
\bea
&&\sum_{n\geq 1}d_3(n)e\left(-\frac{h r n}{q}\right)
V\left(\frac{n}{X}\right)e(\beta r n)\nonumber\\&=&
\frac{q}{2\pi^{\frac{3}{2}}}\sum_{\pm}\sum_{n|q}\sum_{l\geq 1}\frac{1}{n l}
\sum_{n_1|n}\sum_{n_2|\frac{n}{n_1}}\sigma_{0,0}\left(\frac{n}{n_1n_2},l\right)
S\left(-\overline{hr},\pm l;\frac{q}{n}\right)
\Omega_{\beta}^{\pm}\left(\frac{n^2l}{q^3}\right)\nonumber\\
&&+\frac{1}{2q^2}\widetilde{\omega}(1)\sum_{n|q}n\tau(n)P_2(n,q)
\mu\left(\frac{q}{n}\right)
+\frac{1}{2q^2}\widetilde{\omega}'(1)\sum_{n|q}n\tau(n)P_1(n,q)
\mu\left(\frac{q}{n}\right)
\nonumber\\
&&+\frac{1}{4q^2}\widetilde{\omega}''(1)\sum_{n|q}n\tau(n)
\mu\left(\frac{q}{n}\right),
\eea
where $\Omega_{\beta}^{\pm}\left(y\right)$ are defined in (2.9)-(2.10) with
$\omega(y)=V\left(y/X\right)
e(\beta r y)$. As in the Case (i) the first term in (3.11)
is essentially
supported on $n^2l\ll Q^3X^{\epsilon}/X$ and the contribution from the first term
of (3.11) can be bounded similarly as that in
the Case (i), which is at most $\frac{(rX)^{\epsilon}Q^{5/2}}{H}$ with
$Q=(rH)^{2/7}X^{3/7}$.
For the remaining terms in (3.11), we have trivially
\bna
\widetilde{\omega}^{(j)}(1)=\int_0^{\infty}\omega(u)(\log u)^j\mathrm{d}u
\ll_j X(\log X)^j,
\ena
and
they contribute
(3.6) by
\bna
\frac{1}{L}\sum_{\pm}\sum_{1\leq |h|\leq \frac{Q^{1+\epsilon}}{H}}
\sum_{q\in \mathfrak{Q}}\frac{1}{q}
\sum_{|m|\ll Q^2(rX)^{\epsilon}/(rX)}
|a_g(\mp m)|(rX)^{1+\epsilon}\frac{X^{1+\epsilon}}{q^2}\sum_{n|q}n\tau(n)\log^2(nq)
\ll (rX)^{\epsilon}\frac{X}{H}
\ena
which is $\ll \frac{(rX)^{\epsilon}Q^{5/2}}{H}$ for
$Q=(rH)^{2/7}X^{3/7}$.
This finishes the proof of Theorem 1.
\medskip

\section{Proof of Theorem 2}
\setcounter{equation}{0}
\medskip

By dyadic subdivisions we only need to estimate
\bna
\mathcal {T}^{\sharp}(H,Y):=\frac{1}{H}\sum_{h\geq 1}W\left(\frac{h}{H}\right)
\sum_{Y<n\leq 2Y}d_3(n)a_g(rn+h),
\ena
where $Y=2^{-\ell}X$, $1\leq \ell\ll \log X$, $\ell\in \mathbb{Z}$.
Note that
\bna
\mathcal {T}^{\sharp}(H,Y)=\frac{1}{H}\sum_{h\geq 1}W\left(\frac{h}{H}\right)
\sum_{Y<n\leq 2Y}d_3(n)\sum_{m\geq 1}a_g(m)
\phi\left(\frac{m}{rY+h}\right)\delta\left(rn+h-m\right),
\ena
where $\phi(y)$ is as in Theorem 1, i.e., a smooth function compactly supported in $[1/2,5/2]$,
equals 1 on $[1,2]$ and satisfies $\phi^{(j)}(y)\ll_j 1$. Taking $\mathcal {P}=\sqrt{4rY+2H}$ and
applying Lemma 1, we have
\bea
\mathcal {T}^{\sharp}(H,Y)&=&\frac{c_\mathcal {P}}{H\mathcal {P}^2}
\sum_{q\leq \mathcal {P}}\;\sideset{}{^*}\sum_{c \bmod q}\,
\sum_{Y<n\leq 2Y}d_3(n)e\left(\frac{c rn}{q}\right)
\sum_{m\geq 1}a_g(m)e\left(-\frac{cm}{q}\right)\nonumber\\
&&\sum_{h\geq 1}e\left(\frac{ch}{q}\right)W\left(\frac{h}{H}\right)\phi\left(\frac{m}{rY+h}\right)
h\left(\frac{q}{\mathcal {P}},\frac{rn+h-m}{\mathcal {P}^2}\right).
\eea
Then for $H>(rX)^{\frac{1}{2}+\epsilon}$, the proof of Theorem 2 is the similar as
that of Theorem 1 by applying Poisson summation to the $h$-sum in (4.1) and using (3.2).

For $H\leq(rX)^{\frac{1}{2}+\epsilon}$, we let
\bna
\mathcal {T}(H,Y):=\frac{1}{H}\sum_{h\geq 1}W\left(\frac{h}{H}\right)
\sum_{n\geq 1}d_3(n)a_g(rn+h)U\left(\frac{n}{Y}\right),
\ena
where $U(y)$ is a smooth function compactly supported in $[1,2]$, which equals 1
on $[1+\Delta^{-1},2-\Delta^{-1}]$ ($\Delta>1$ is a parameter to be chosen
optimally later) and satisfies $U^{(j)}(y)\ll_j \Delta^j$.
Assume that $a_g(n)\ll n^{\theta+\epsilon}$. Then we have
\bea
\mathcal {T}^{\sharp}(H,Y)=\mathcal {T}(H,Y)
+O_{g,\epsilon}\left(X\Delta^{-1}(r X)^{\theta+\epsilon}\right).
\eea
Note that
\bna
\mathcal {T}(H,Y)=\frac{1}{H}\sum_{h\geq 1}W\left(\frac{h}{H}\right)
\sum_{n\geq 1}d_3(n)U\left(\frac{n}{Y}\right)\sum_{m\geq 1}a_g(m)
\phi\left(\frac{m}{rY+h}\right)\int_0^1e\left((rn+h-m)\alpha\right)\mathrm{d}\alpha.
\ena
As in the proof of Theorem 1 we apply Lemma 2 to approximate $\mathcal {T}(H,Y)$ by
\bna
\widetilde{\mathcal {T}}(H,Y)&=&\frac{1}{H}\sum_{h\geq 1}W\left(\frac{h}{H}\right)\int_0^1
\widetilde{I}_{\mathfrak{Q},\eta}(\alpha)e(\alpha h)
\sum_{n\geq 1}d_3(n)e(\alpha r n)U\left(\frac{n}{Y}\right)\\
&&\sum_{m\geq 1}a_g(m)
e(-\alpha m)\phi\left(\frac{m}{rY+h}\right)\mathrm{d}\alpha,
\ena
where $\widetilde{I}_{\mathfrak{Q},\eta}(\alpha)(x)$ is defined in (2.3) with
\bna
\mathfrak{Q}=\{q: q\in [\mathcal {Q}/2,\mathcal {Q}] \ \text{is prime and}\, (q,r)=1\}.
\ena
Take $\eta=(rY+H)^{-1}$. Then by Cauchy's inequality and (2.4),
\bea
\mathcal {T}(H,Y)=\widetilde{\mathcal {T}}(H,Y)
+O\left((r X)^{\epsilon}\frac{rX^{\frac{3}{2}}}{\mathcal {Q}}\right).
\eea
In the following we estimate $\widetilde{\mathcal {T}}(H,Y)$.
We have
\bna
\widetilde{\mathcal {T}}(H,Y)=\frac{1}{2\eta}\int_{-\eta}^{\eta}\widetilde{\mathcal {T}}_{\beta}(H,Y)\mathrm{d}\beta,
\ena
where
\bea
\widetilde{\mathcal {T}}_{\beta}(H,Y)&=&
\frac{1}{HL}\sum_{q\in \mathfrak{Q}}\;\sideset{}{^*}\sum_{c \bmod q}
\sum_{n\geq 1}d_3(n)e\left(\frac{c r n}{q}\right)
U\left(\frac{n}{Y}\right)e\left(\beta r n\right)\sum_{m\geq 1}a_g(m)
e\left(-\frac{c m}{q}\right)
e\left(-\beta m\right)\nonumber\\
&&\sum_{h\geq 1}e\left(\frac{c h}{q}\right)W\left(\frac{h}{H}\right)
\phi\left(\frac{m}{rY+h}\right)e\left(\beta h\right).
\eea
As in Theorem 1, we apply Poisson summation to the $h$-sum to get
\bea
\mbox{$h$-sum}=H\sum_{1\leq |h|\leq \frac{\mathcal {Q}^{1+\epsilon}}{H}\atop h\equiv -c\bmod q}I_{\beta}(h,m,q)
+O((r X)^{-A})
\eea
for any $A>0$. Plugging (4.5) into (4.4), we have
\bna
&&\frac{1}{L}\sum_{1\leq |h|\leq \frac{Q^{1+\epsilon}}{H}}
\sum_{q\in \mathfrak{Q}\atop (q,h)=1}
\sum_{n\geq 1}d_3(n)e\left(-\frac{h r n}{q}\right)
U\left(\frac{n}{Y}\right)e\left(\beta rn\right)\sum_{m\geq 1}a_g(m)
e\left(\frac{h m}{q}\right)
e\left(-\beta m\right)I_{\beta}(h,m,q)\\
&:=&\int_{\mathbb{R}}W\left(x\right)
e\left(\beta H x\right)
\widetilde{\mathcal {T}}_{\beta,x}(H,Y)\mathrm{d}x,
\ena
where
\bea
\widetilde{\mathcal {T}}_{\beta,x}(H,Y)&=&
\frac{1}{L}\sum_{1\leq |h|\leq \frac{Q^{1+\epsilon}}{H}}
\sum_{q\in \mathfrak{Q}\atop (q,h)=1}e\left(-\frac{h H x}{q}\right)
\sum_{n\geq 1}d_3(n)e\left(-\frac{h rn}{q}\right)
U\left(\frac{n}{Y}\right)e\left(\beta r n\right)\nonumber\\
&&\sum_{m\geq 1}a_g(m)
e\left(\frac{h m}{q}\right)\phi\left(\frac{m}{rY+Hx}\right)
e\left(-\beta m\right).
\eea
Applying the $GL(2)$ Voronoi formula in Lemma 2 to
the $m$-sum in (4.6) we get
\bea
\widetilde{\mathcal {T}}_{\beta,x}(H,Y)&=&
\frac{1}{L}\sum_{\pm}\sum_{1\leq |h|\leq \frac{\mathcal {Q}^{1+\epsilon}}{H}}
\sum_{q\in \mathfrak{Q}\atop (q,h)=1}\frac{1}{q}e\left(-\frac{h H x}{q}\right)
\sum_{|m|\ll \mathcal {Q}^2(rY+H)^{\epsilon}/(rY+H)}
a_g(\mp m)e\left(\pm \frac{\overline{h} m}{q}\right)\nonumber\\
&&\Psi_{\beta,x}^{\pm}\left(\frac{m}{q^2}\right)
\sum_{n\geq 1}d_3(n)e\left(-\frac{h r n}{q}\right)
U\left(\frac{n}{Y}\right)e\left(\beta r n\right)
+O_{g,\epsilon}(1).
\eea
where $\Psi_{\beta,x}^{\pm}\left(y\right)$ are defined in (2.5)-(2.6) with
$\psi(y)=\phi\left(\frac{y}{rY+Hx}\right)
e\left(-\beta y\right)$.
and satisfy $\Psi_{\beta,x}^{\pm}\left(\frac{m}{q^2}\right)\ll (rY+H)^{1+\epsilon}$.

Next we apply the Voronoi formula
for $d_3(n)$ in Lemma 5 to the $n$-sum
in (4.7) to get (3.11) with  $V\left(\frac{n}{X}\right)$
replaced by $U\left(\frac{n}{Y}\right)$ and
$\omega(y)=U\left(\frac{y}{Y}\right)e(\beta ry)$.
By Weil's bound for Kloosterman sums, the contribution from the last three terms
in (3.11) to (4.7) is at most
\bea
&&\frac{1}{L}\sum_{\pm}\sum_{1\leq |h|\leq \frac{\mathcal {Q}^{1+\epsilon}}{H}}
\sum_{q\in \mathfrak{Q}\atop (q,h)=1}\frac{1}{q}
\sum_{|m|\ll \mathcal {Q}^2(rY+H)^{\epsilon}/(rY+H)}
|a_g(\mp m)|\frac{(rY+H)^{1+\epsilon}Y}{q^2}\sum_{n|q}n\tau(n)
(\log nq)^2\nonumber\\
&\ll_{g,\epsilon}&(rX)^{\epsilon}\frac{X}{H}.
\eea

For the first term in (3.11), we note that
$
\frac{d^j}{dy^j}\left\{U\left(\frac{y}{Y}\right)
e(\beta ry)\right\}\ll_{g,\epsilon}
\left(\frac{\Delta}{Y}\right)^{j}
$
for any $j\geq 0$. By Lemma 6, the contribution from
$n^2l\gg q^3\Delta^3(qY)^{\epsilon}/Y$ is negligible.
For $n^2l\ll q^3\Delta^3Y^{\epsilon}/Y$, we shift the contour of integration
in (2.9) to $\sigma=-1/2-\epsilon$ with $\epsilon>0$ to get
\bea
\Omega_{\beta}^{\pm}\left(\frac{n^2l}{q^3}\right)\ll_{\epsilon}\Delta
\left(\frac{Yn^2l}{q^3}\right)^{1/2+\epsilon}.
\eea
By (4.9) and Weil's bound for Kloosterman sums, one sees that the
first term in (3.11) contributes $\widetilde{\mathcal {T}}_{\beta,x}(H,Y)$ in (4.7) by
\bea
&&\frac{1}{L}\sum_{\pm}\sum_{1\leq |h|\leq \frac{\mathcal {Q}^{1+\epsilon}}{H}}
\sum_{q\in \mathfrak{Q}\atop (q,h)=1}
\sum_{|m|\ll \mathcal {Q}^2(rY+H)^{\epsilon}/(rY+H)}
|a_g(\mp m)|(rY+H)^{1+\epsilon}\nonumber\\
&&\sum_{n|q}\sum_{l\ll q^3\Delta^3Y^{\epsilon}/(n^2Y)}\frac{1}{nl}
\sum_{n_1|n}\sum_{n_2|\frac{n}{n_1}}\sigma_{0,0}\left(\frac{n}{n_1n_2},l\right)
\left(\frac{q}{n}\right)^{1/2}
\left(\frac{Yn^2l}{q^3}\right)^{1/2+\epsilon}\nonumber\\
&\ll_{g,\epsilon}&
\frac{(rY+H)^{\epsilon}\mathcal {Q}}{H}
\sum_{\mathcal {Q}/2\leq q\leq \mathcal {Q} }\sum_{n|q}\,\sum_{l\ll q^3\Delta^3(qY)^{\epsilon}/(n^2Y)}
\frac{d_3(l)d_3(n)}{nl}\left(\frac{q}{n}\right)^{1/2}
\left(\frac{Yn^2l}{q^3}\right)^{1/2}\nonumber\\
&\ll_{g,\epsilon}&
\frac{(rY+H)^{\epsilon}\mathcal {Q}Y^{1/2}\Delta}{H}
\sum_{\mathcal {Q}/2\leq q\leq \mathcal {Q} }q^{-1}\sum_{n|q}n^{-1/2}
\sum_{l\ll q^3\Delta^3(qY)^{\epsilon}/(n^2Y)}
l^{-1/2}\nonumber\\
&\ll_{g,\epsilon}&(rX)^{\epsilon}
\frac{Q^{5/2}\Delta^{5/2}}{H}.
\eea
By (4.2), (4.3), (4.8) and (4.10) we take $\mathcal {Q}=(rH)^{2/7}X^{3/7}\Delta^{-5/7}$ and obtain
\bna
\mathcal {T}^{\sharp}(H,Y)\ll_{g,\epsilon}
\frac{(rH\Delta)^{5/7}X^{15/14}}{H}+
\frac{X(r X)^{\theta+\epsilon}}{\Delta}.
\ena
Then Theorem 2 follows by choosing $\Delta=(HX)^{1/6}(rX)^{7/(12\theta)}(r^2X)^{-5/24}$.

\medskip

\bigskip

\noindent {\bf Acknowledgments} The author would like to thank Yongxiao Lin for useful
discussions and Department of Mathematics, The Ohio State University for
hospitality. This work is supported by
the National Natural Science Foundation of China (Grant No. 11101239),
Young Scholars Program of Shandong University, Weihai (Grant No. 2015WHWLJH04),
the Natural Science Foundation of Shandong Province (Grant No. ZR2016AQ15)
and a scholarship from the China
Scholarship Council.

\bigskip

{
Qingfeng Sun

School of Mathematics and Statistics

Shandong University, Weihai

Weihai, Shandong 264209

China

{\it email: qfsun@sdu.edu.cn}}

\end{document}